\documentclass[3p]{elsarticle}

\usepackage{graphicx}
\usepackage{amssymb}
\usepackage{amsthm}
\usepackage{lineno}
\usepackage[usenames]{color}

\newdefinition{rmk}{Remark}
\newproof{pf}{Proof}
\newproof{pot}{Proof of Theorem \ref{thm2}}
\begin{document}
\begin{frontmatter}
\title{High-performance modeling acoustic and elastic waves using the Parallel Dichotomy Algorithm.}

\author[geo]{Alexey G. Fatyanov}
\ead{fat@nmsf.sscc.ru}

\author[geo,budker,unv]{Andrew V. Terekhov\corref{cor1}}
\ead{andrew.terekhov@mail.ru}

 \newcommand{\edt}[1]{\{\textcolor{red}{#1}\}}

\address[geo]{Institute of
Computational Mathematics and Mathematical Geophysics,
630090,Novosibirsk,Russia}
\address[budker]{Budker Institute of Nuclear Physics, 630090, Novosibirsk,
Russia}
\address[unv]{Novosibirsk State University, 630090, Novosibirsk,
Russia}

\cortext[cor1]{Corresponding author}

\begin{abstract}
A high-performance parallel algorithm is proposed for modeling the
propagation of acoustic and elastic waves in inhomogeneous media.
An initial boundary-value problem is replaced by a series of
boundary-value problems for a constant elliptic operator and
different right-hand sides via the integral Laguerre transform. It
is proposed to solve difference equations by the conjugate
gradient method for acoustic equations and by the GMRES$(k)$
method for modeling elastic waves. A preconditioning operator was
the Laplace operator that is inverted using the variable
separation method. The novelty of the proposed algorithm is using
the Dichotomy Algorithm (Terekhov, 2010), which was designed for
solving a series of tridiagonal systems of linear equations, in
the context of the preconditioning operator inversion. Via
considering analytical solutions, it is shown that modeling wave
processes for long instants of time requires high-resolution
meshes. The proposed parallel fine-mesh algorithm enabled to solve
real application seismic problems in acceptable time and with high
accuracy. By solving model problems, it is demonstrated that the
considered parallel algorithm possesses high performance and
efficiency over a wide range of the number of processors (from $2$
to $8192$).
\end{abstract}

\end{frontmatter}

\section{Introduction}
Steadily growing number of processors opens up new opportunities
for solving complex applied problems, for example, elastodynamic
problems \cite{Kaufman3,Achenbach,White}. In this case, quite
efficient algorithms are numerical-analytical algorithms
\cite{Mikhailenko1999,Flayer2}, where the solution is represented,
via the integral time transformation, as the Fourier series in
terms of some orthonormal system of functions. The expansion
coefficients are determined numerically as a solution of
boundary-value problems for the elliptic type of equation
\cite{Samarskii2001,FEM2,FEM1}.

Many publications
\cite{Ell_prob_4,Ell_prob_5,Ell_prob_1,Ell_prob_2,Ell_prob_3,Ell_prob_6},
are concerned with development and investigation of parallel
numerical elliptic differential operator inversion algorithms.
Nevertheless, this problem remains quite urgent. The explanation
is that the steady growth of the number of processors integrated
within one computer system imposes new demands to scalability of
parallel algorithms. For instance, methods effective for a small
number of processors ($p<32$), e.g., the cyclic reduction
algorithm  \cite{hockney:cyclic,Johnsson:cyclic}, become
ineffective because the communication costs prevail over the
computational ones. This necessitates further development of
parallel numerical algorithms that allow using modern
computational resources with the greatest efficient factor.

Publications
\cite{Samarski_Nikolaev,bernhardt,Concus:Golub,Israeli,Chang,Waveform}
propose different approaches to solving elliptic equations of
second order with inseparable variables, where the iterative
process is reduced to multiple Laplacian inversion. However,
realization of efficient procedure
\cite{Samarski_Nikolaev,Hockney:Particles,Hockney,Pissman}for the
Laplacian inversion requires solving tridiagonal systems of linear
equations, which, on a multiprocessor system, is a nontrivial
problem. This difficulty can be overcome by using the Dichotomy
Algorithm \cite{Terekhov:dichotomy}, which was developed for
inversion of one and same tridiagonal matrix for many right-hand
sides. The Dichotomy Algorithm was chosen because for this class
of problems it ensures almost linear dependence of the speedup
coefficient in a wide range of the number of processors. In terms
of accuracy, the number of arithmetic operations, and the number
of communications, the Dichotomy Algorithm is practically
equivalent to the cyclic reduction method
\cite{hockney:cyclic,Snonkwiler:Parallel}. However, with
comparable levels of transferred data, the real time of
interprocessor communications for the Dichotomy Algorithms is much
less. The explanation is that the main communication operation of
the Dichotomy Algorithm, that is, all-reduce-to-one(+), possesses
the associative property, which allows reducing the time of
interprocessor communications due to their optimization
\cite{Tuning:2,Collective:1}. In the present paper, taking into
account the high efficiency of the Dichotomy Algorithm, we will
consider the possibility of using it within the scope of
numerical-analytical approach for modeling the propagation of
acoustic and elastic waves.

The peculiarity of the Dichotomy Algorithm is that it was designed
for solving problems with the same tridiagonal matrix and
different right-hand sides. Choosing the integral transformation,
we considered the fact that prior to solving tridiagonal systems,
it is necessary to perform preparatory calculations with a volume
$O(N)$, where  $N$ is the dimension of the system of equations.
Really, after applying the time Fourier transform to the acoustic
equation we obtain the boundary-value problem for the Helmholtz
equation
\begin{equation} \label{Helmholtz} \Delta
u_n+k^2_nu_n=f_n,\quad n=1,2,...\;.
\end{equation}
In this case, the dependence of the differential operator on the
number of calculated harmonic prevents effective usage of the
Dichotomy Algorithm because only one right-hand side will
correspond to the same matrix. The exception is the case of
Toeplitz tridiagonal matrices \cite{Terekhov2} for which the
volume of preparatory computations is  $O(N/p+\log_2 p)$ rather
than $O(N)$, where $p$ is the number of processors. Thus, for
solution of problem (\ref{Helmholtz}) in the Cartesian coordinate
system, the Dichotomy Algorithm can be applied, e.g., in the
context of the variable separation method that requires inversion
of Toeplitz (quasi-Toeplitz) matrices.

In the present work, we consider media of 2.5D geometry. In this
case, in the cylindrical coordinate system, for the Laplace
operator inversion, it is necessary to solve tridiagonal SLAEs of
the general form. For this case, we considered the Laguerre
transform\cite{Mikhailenko1999}, after applying it to the acoustic
equation, it is required to invert one and the same differential
operator for all right-hand sides

\begin{equation}
\label{Laguerre_exmpl} \Delta
u_n-\lambda^2u_n=f_n+\sum_{i=1}^{n-1}\alpha_{n,i}u_i,\quad
n=1,2,... \quad\quad \alpha_{n,i},\lambda \in R.
\end{equation}

The fact that the preparatory computations in the context of the
Dichotomy Algorithm are performed once for all right-hand sides
allows one to neglect preparatory expenses. Thus, it becomes
possible to use the Dichotomy Algorithm for solving problem
(\ref{Laguerre_exmpl}) by methods demanding inversion of general
tridiagonal matrices.

In the present paper, using the Laguerre transform and the
Dichotomy Algorithm we considered the high-performance parallel
algorithm for modeling acoustic and elastic waves in 2.5D media.
At present, applied geophysics problems have to be solved for
steadily increasing times and recording systems. On the other
hand, improvement of practical observing systems necessitates
increasing the calculation accuracy. In the paper, by considering
analytical solutions we have shown that modeling wave processes
for longer instants of time requires higher resolution meshes. We
illustrated the possibility of effective using thousands of
processors within one calculation. This enabled practical
real-time and high-accuracy computing on current computers.

\section{Problem Statement and Solution Algorithm}
\subsection{Acoustic Equation}
In the cylindrical coordinate system  $(r,z)$, in the half-space
$z\geq 0$ we will consider the problem of modeling the propagation
of acoustic waves from a point source

\begin{equation}
\label{acoustic_problem}
\begin{array}{llr}
\displaystyle {\rho({\bf x} )}\frac{\partial^2 u}{\partial
t^2}({\bf x},t)=\nabla \left[\kappa({\bf x})\, \nabla u({\bf
x},t)\right]+\frac{1}{2\pi}\frac{\delta({\bf x-x_0})}{r}f(t),&
t>0,\quad {\bf x}=(r,z).
\end{array}
\label{wave-eq}
\end{equation}

Suppose that problem (\ref{acoustic_problem}) is solved with
homogeneous initial conditions
\begin{equation}
\begin{array}{llr}
\displaystyle \left. u\right|_{t=0}=\left.\frac{\partial
u}{\partial t}\right|_{t=0}=0.
\end{array}
\label{initcond}
\end{equation}
Assume that at $z=0$ the surface is free, and the auxiliary
boundaries are entered along the coordinates $r$ and $z$
\begin{equation}
\begin{array}{llr}
\displaystyle \left.\frac{\partial u}{\partial
z}\right|_{z=0,l_2}=\left.u\right|_{r=l_1}=0.
\end{array}
\label{boundary_cond}
\end{equation}

The boundaries $r=l_1$ and $z=l_2$ are chosen such that waves
reflected from them do not arise for the calculated instant of
time. In addition we demand that
\begin{equation} \left.\frac{\partial
u}{\partial r}\right|_{r=0}=0. \label{symm_cond}
\end{equation}

Let us represent for the solution of problem
(\ref{acoustic_problem})-(\ref{symm_cond}) as the Fourier-Laguerre
series\cite{Mikhailenko1999}

\begin{equation}
R_m(\mathbf{x})=\int_{0}^{\infty}u(\mathbf{x},t)(ht)^{-\frac{\alpha}{2}}l_{m}^{\alpha}(ht)dt
\label{series_lag}
\end{equation}
with the inversion formulas
\begin{equation}
u(\mathbf{x},t)=(ht)^{\frac{\alpha}{2}}\sum_{m=0}^{\infty}R_m({\bf
x})l^{\alpha}_m(ht) , \label{series_lag}
\end{equation}

where $l^{\alpha}_m(ht)$ are the orthonormal Laguerre functions
\cite{abramowitz+stegun}, which are represented via classical
Laguerre polynomials as follows
$$
l^{\alpha}_m(ht)=\sqrt{\frac{hm!}{(m+\alpha)!}}(ht)^{\frac{\alpha}{2}}e^{-\frac{ht}{2}}L^{\alpha}_m(ht).
$$

Here, $m$ is Laguerre polynomial degree and $h$ is the
transformation parameter. The necessary and sufficient parameter
for satisfying the initial data is $\alpha \ge 2$ ($\alpha$ is the
order of Laguerre functions).

As a result, the initial boundary-value problem
(\ref{acoustic_problem})--(\ref{symm_cond}) is reduced to the
boundary-value problems in the spectral domain

\begin{equation}
\left\{
\begin{array}{l}
\label{laguerre_h} \displaystyle \nabla \left[\kappa({\bf x})\,
\nabla R_m({\bf x})\right]-\rho({\bf x})\frac{h^2}{4}R_m({\bf
x})=-\frac{1}{2\pi}\frac{\delta({\bf
x-x_0})}{r}f_m+\rho({\bf x})h^2\sqrt{\frac{m!}{(m+\alpha)!}}\sum_{k=0}^{m-1}(m-k)\sqrt{\frac{(k+\alpha)!}{k!}}R_k({\bf x}),\\\\
\displaystyle\left. \frac{\partial R_m}{\partial
r}\right|_{r=0}=\left.\frac{\partial R_m}{\partial
z}\right|_{z=0,l_2}=\left.R_m\right|_{r=l_1}=0,
\end{array}\right.
\end{equation}
where
$f_m=\int_0^{\infty}f(t)(ht)^{-\frac{\alpha}{2}}l^{\alpha}_{m}(ht)dt$.

This method can be considered as an analog of the
spectral-difference method, based on the Fourier transform
\cite{FEM1}, but in this case, but the role of "frequency"\
belongs to the parameter $m$ that determines the degree of the
polynomials. Contrary to the Fourier method, the harmonic
separation parameter is present only in the right-hand side.
\subsection{Elastic Medium}
To describe the propagation of elastic waves in a inhomogeneous
half-space, we will consider the equations of motion in the
cylindrical coordinate system \cite{White}
\begin{equation}
\label{elastic_problem}
\begin{array}{l}
\displaystyle {\rho}\frac{\partial \mathbf{W}^2}{\partial t^2}=
%\sigma_{ij,j}+\mathbf{F}_i(\textbf{x},t)
%\\\\
\left(\lambda+\mu\right)\nabla\left(\nabla\cdot
\mathbf{W}\right)+\mu\nabla^2 \mathbf{W}+\nabla\lambda\left(\nabla
\cdot\mathbf{W}\right)+\nabla\mu\times\left(\nabla\times\mathbf{W}\right)+2\left(\nabla\mu\cdot\nabla\right)\mathbf{W}+\rho\mathbf{F}.
\label{boundary_elastic}
\end{array}
\end{equation}
Here, $\mathbf{W}$ is the displacement vector, $\lambda>0$ and
$\mu>0$ are Lame coefficients, $\mathbf{F}$ is the force vector
describing the action of space-localized axially symmetric source.

Let us consider the case of the cylindrical coordinate system
(2.5D), where $\mathbf{W}=(u_r,u_z)^{\mathrm{T}}$,
$\mathbf{F}=\left(F_r,F_z\right)^{\mathrm{T}}$,
$\lambda=\lambda(r,z),\mu=\mu(r,z)$ and $\rho=\rho(r,z)$. Assume
that at $z=0$ , the surface is free, with auxiliary boundaries
along the coordinates $r$ and $z$, as in the case of the acoustic
equation. Problem (\ref{elastic_problem}) is solved with
homogeneous initial conditions.

Represent the solution of problem  (\ref{elastic_problem}) as the
Fourier-Laguerre series

\begin{equation}
\label{series_lag2}
u_r(\mathbf{x},t)=(ht)^{\frac{\alpha}{2}}\sum_{m=0}^{\infty}Q_m({\bf
x})l^{\alpha}_m(ht),\quad
u_z(\mathbf{x},t)=(ht)^{\frac{\alpha}{2}}\sum_{m=0}^{\infty}U_m({\bf
x})l^{\alpha}_m(ht) .
\end{equation}

As a result, defining the expansion coefficients $Q_m$ and $U_m$
necessitates solving a number of problems of the form

\begin{equation}
\label{elastic_problemQn}
\left\{
\begin{array}{l}
\displaystyle \frac{\partial}{\partial
r}\left[\left(2\mu+\lambda\right)\frac{\partial Q_m}{\partial r}
+\lambda\left( \frac{
\partial U_m}{\partial z}+\frac{Q_m}{r}\right)\right]+\frac{\partial}{\partial
z}\left[\mu \left( \frac{\partial Q_m}{\partial z}+\frac{\partial
U_m}{\partial r}\right)\right]+\frac{2\mu}{r}\left(\frac{\partial
Q_m}{\partial r}-\frac{Q_m}{r}\right)-\rho\frac{h^2}{4}Q_m=\\\\
\displaystyle =-\rho F_rf_m+\rho
h^2\sqrt{\frac{m!}{(m+\alpha)!}}\sum_{k=0}^{m-1}(m-k)\sqrt{\frac{(k+\alpha)!}{k!}}Q_k,
\\\\
\displaystyle\frac{1}{r}\frac{\partial}{\partial r}\left[r\mu
\left(\frac{\partial Q_m}{\partial z}+\frac{\partial U_m}{\partial
r}\right) \right]+\frac{\partial}{\partial z}\left[\left(\lambda+2
\mu\right) \frac{\partial U_m}{\partial z}+\lambda\left(
\frac{\partial U_m}{\partial r}
+\frac{Q_m}{r}\right)\right]-\rho\frac{h^2}{4}U_m=\\\\
\displaystyle = -\rho F_zf_m+\rho
h^2\sqrt{\frac{m!}{(m+\alpha)!}}\sum_{k=0}^{m-1}(m-k)\sqrt{\frac{(k+\alpha)!}{k!}}U_k,
\end{array}\right.
\end{equation}

where the boundary conditions on the free surface take the
form\cite{Kaufman3,Achenbach,White}

\begin{equation}
\tilde{\tau}_{rz}=\left\{\left .\frac{\partial U_m}{\partial
z}+\frac{\partial Q_m}{\partial r}\right\}\right|_{z=0}=0,
\end{equation}

\begin{equation}
\tilde{\sigma}_{zz}=\left\{\left . \lambda\left(\frac{\partial
U_m}{\partial
r}+\frac{U_m}{r}\right)+\left(\lambda+2\mu\right)\frac{\partial
Q_m}{\partial z}\right\}\right|_{z=0}=0.
\end{equation}

\subsection{Approximation of equations}
 On a rectangular mesh $ \label{grid2}
\bar{\omega}=\bar{\omega}_r\times\bar{\omega}_z=\omega\bigcup\gamma$
,where
$$
\begin{array}{l} \bar \omega_r=\left\{r_i=(i-0.5)h_r,\;
i=1,...,N_r,\; h_r=l_1/(N_r-0.5)\right\},\\\\
\bar \omega_z=\left\{z_k=(k-0.5)h_z,\; k=1,...,N_z,\;
h_z=l_2/(N_z-0.5)
\right\},\\\\
\omega=\bar\omega\bigcap G,\quad \gamma=\bar\omega\bigcap\Gamma,
\end{array}
$$

conform problem  (\ref{laguerre_h}) with the difference problem

\begin{equation}
\label{acoustic_problem_A}
 Ay_m=f,\; m=1,2,..., \quad A:H\longrightarrow H,
\end{equation}

where the difference operator $A=A^{*}>0$ is given by a scheme of
the second order of accuracy \cite{Samarskii2001,Strikwerda2004}

\begin{equation}
\begin{array}{ll}
\displaystyle\label{poisson_diff_2}
\left(\Lambda_r+\Lambda_z\right)y_m-w(x)y_m=-\phi(x),\quad x \in
\bar{\omega},
\end{array}
\end{equation}

\begin{equation}
\label{scheme_op_r_z}
\begin{array}{lr}
\Lambda_ry=\left\{
\begin{array}{ll}
\displaystyle \frac{1}{h_r}a_1y_{r},& i=1\\\\
\displaystyle \left(a_1y_{\bar{r}}\right)_{r},&1\leq i\leq N_1-1
\end{array}\right., &
\Lambda_zy=\left\{
\begin{array}{ll}
\displaystyle \frac{1}{h_z}a_2y_{z},& k=1\\\\
\displaystyle \left(a_2y_{\bar{z}}\right)_{z},&1\leq k\leq
N_2-1\\\\
\displaystyle \frac{1}{h_z}a_2y_{\bar{z}},&k=N_2
\end{array},\right.
\end{array}
\end{equation}

\begin{equation}
\label{diff_approx1}
\begin{array}{l}
 \displaystyle
a_1(i,k)=\bar{r}_i\kappa\left(\bar{r}_i,z_k\right),\;
a_2(i,k)=r_i\kappa\left(r_i,\bar{z}_k\right),\;
w(i,k)=\rho(r_i,z_k)\frac{h^2}{4}r_i,\\\\
\phi(i,j)=\displaystyle
-\frac{1}{2\pi}\frac{\delta(i-i_0,j-j_0)}{r_i}f_m+\rho({\bf
x})h^2\sqrt{\frac{m!}{(m+\alpha)!}}\sum_{k=0}^{m-1}(m-k)\sqrt{\frac{(k+\alpha)!}{k!}}y_k({\bf
x}) ,\quad {\bf x}_0=(i_0h_r,j_0h_z).
%\left\{\begin{array}{cl}\displaystyle\frac{f_m}{2\pi
%h_r
%hz},& if\; {\bf x}_0=(ih_r,j h_z)
\\\\
%0,& else
%\end{array}\right.%\\
\end{array}
\end{equation}

where $\bar{r}_i=r_i+0.5h_r,\; \bar{z}_k=z_k+0.5h_z$; $\;y_{\bar
r},\,y_{\bar z}$ and $y_r,\,y_z$  are the "backward"\  and
"forward"\ difference relationships with respect to $z$ and $r$
\cite{Samarski_Nikolaev,Samarskii2001}. The boundary condition on
the side $r=l_1$ is approximated exactly $ y_{N_1,k}=0, \quad
k=1,...,N_2$. For solving problem (\ref{poisson_diff_2}) , we use
the conjugate gradient method \cite{Saad}.

Conform the problem (\ref{elastic_problem}) on a mesh
$\bar{\omega}$ with the difference problem

\begin{equation}
Cy_m=f,\; m=1,2,..., \quad C:H\rightarrow H.
\label{elastic_problem_slau}
\end{equation}

A number of works \cite{Approx:Mech1,Approx:Mech2} describe the
problem of constructing the discrete analog of problem
(\ref{elastic_problemQn}). For this reason, in the present paper
we performed approximation by the finite-volume method with second
order of accuracy. For solving problem
(\ref{elastic_problem_slau}), by virtue of the non self-adjoint
difference operator $C$, we will use the GMRES$(k)$ method, where
$k$ is the restart parameter \cite{Saad}. Note, when using the
Laguerre transformation, the difference operator is always
positive-definite. This guarantees convergence of the GMRES$(k)$
method for any $k\geq1$\cite{Saad}.

\subsection{Preconditioning}
By choosing a preconditioning procedure, one can affect
substantially the convergence of iterative algorithms of solving a
system of linear equations and, as a result, the elapsed time.
Besides standard requirements \cite{Samarski_Nikolaev,Saad} upon a
preconditioning operator, such as

\begin{itemize}
\item energy equivalence of operator $B$ to operator $A$ in the
sense of inequalities \footnote{For the non self-adjoint case,
see, e.g., \cite{Samarski_Nikolaev,Saad,Dongarra:NumericalLA}.}
\begin{equation}
    \label{energy_rel}
    \gamma_1 \left(Bu,u\right) \leq \left(Au,u\right) \leq
    \gamma_2 \left(Bu,u\right);  \quad \quad
    0<\gamma_1\leq\gamma_2,\quad A=A^*>0,\,B=B^*>0,
    \end{equation}
    $$\gamma_1=\min_{x\neq
     0}\frac{\left(Ax,x\right)}{\left(Bx,x\right)},\quad
     \gamma_2=\max_{x\neq 0}\frac{\left(Ax,x\right)}{\left(Bx,x\right)};$$
    \item  operation of inversion of operator $B$ must be less time-consuming than for operator $A$,
\end{itemize}

we require efficiency of a procedure preconditioning of operator
inversion on a multiprocessor computer system. Since not all
preconditioning procedures can be efficiently implemented with the
use of hundreds of processors (e.g., ILU expansion \cite{Saad}),
the latter requirement drastically limits the class of possible
preconditioners.

In \cite{Terekhov:dichotomy}, based on the Dichotomy Algorithm,
the author proposed a high-performance parallel implementation of
the variable separation method
\cite{Samarski_Nikolaev,Hockney,FFT} for the Laplace operator
inversion. The use of the Dichotomy Algorithm for solving
tridiagonal systems of linear equations ensures linear dependence
of the speedup coefficient on the number of processors. Thus, for
problem (\ref{acoustic_problem_A}), following works
\cite{Samarski_Nikolaev,bernhardt,Concus:Golub,Israeli,Chang,Waveform},
as the preconditioner operator we will consider\footnote{Introduce
the notation $\widetilde{f}=\frac{1}{2}\left(\min_{x \in G }f({\bf
x})+\max_{x \in G }f({\bf x})\right).$}

\begin{equation}
\label{precond1} B\equiv \Lambda_r+\Lambda_z-d
\end{equation}
with the coefficients
$$a_1(i,k)=\bar{r}_i\widetilde{{\kappa}},\quad
a_2(i,k)=r_i\widetilde{\kappa},\quad \displaystyle
d(i,k)=r_i\frac{h^2}{4}\widetilde{\rho}.$$

For problem (\ref{elastic_problem_slau}), the preconditioner is
given as

\begin{equation}
\label{precond2} K\equiv\left[
\begin{array}{cc}
B_1&0 \\\\
0&B_2
\end{array}
\right],
\end{equation}
where $\displaystyle B_1 \equiv \Lambda_r+\Lambda_z-d$ ñ with the
coefficients
$$a_1(i,k)=\bar{r}_i\widetilde{(\lambda+2\mu)},\quad
a_2(i,k)=r_i\widetilde{\mu},\quad \displaystyle
d(i,k)=r_i\frac{h^2}{4}\widetilde{\rho^{}}+\frac{\widetilde{(\lambda+2\mu)}}{r_i}$$

and $B_2 \equiv \Lambda_r+\Lambda_z-d$ with the coefficients
$$a_1(i,k)=\bar{r}_i\widetilde{\mu},\quad
a_2(i,k)=r_i\widetilde{(\lambda+2\mu)},\quad \displaystyle
d(i,k)=r_i\frac{h^2}{4}\widetilde{\rho^{}}$$

By virtue of the assumption that the contrast of the medium is
moderate and the use of a supercomputer implies a great number of
mesh nodes, this class of preconditioners enables a good
convergence rate. Moreover, the sought solution will be achieved
in the number of iterations, which does not practically depend on
the number of mesh nodes \cite{Samarski_Nikolaev}.

Since in the proposed scheme of solution of to problem the main
computational and communication costs fall on the preconditioner
inversion, the algorithm efficiency, as a whole, is determined by
performance of the parallel procedure of solution of the problems
$B_{\alpha}y=\phi$.

\section{Numerical Experiments}
\subsection{Parallel Performance}
For estimating the performance of the proposed algorithm, using
Fortran-90 and the MPI paradigm, we implemented numerical
procedures for solving problems (\ref{acoustic_problem}) and
(\ref{elastic_problem}). The Fast Fourier transform, which is
necessary for the preconditioner inversion, was done using FFTW
library \cite{FFTW}; the tridiagonal systems of linear equations
were solved using the Dichotomy Algorithm
\cite{Terekhov:dichotomy,Terekhov2}. Calculations were performed
on MBC-100k supercomputer (from the Interdepartment Supercomputer
Center of the Russian Academy of Sciences) and on NKS-30t
supercomputer (from the Siberian Supercomputer Center of the
Siberian Branch of the Russian Academy of Sciences). The computer
are based on Intel Xeon four-core processors operating at $3$ GHz
and connected via the Infiniband communication medium.

\begin{figure}[!h]
\begin{center}
\includegraphics[width=0.5\textwidth]{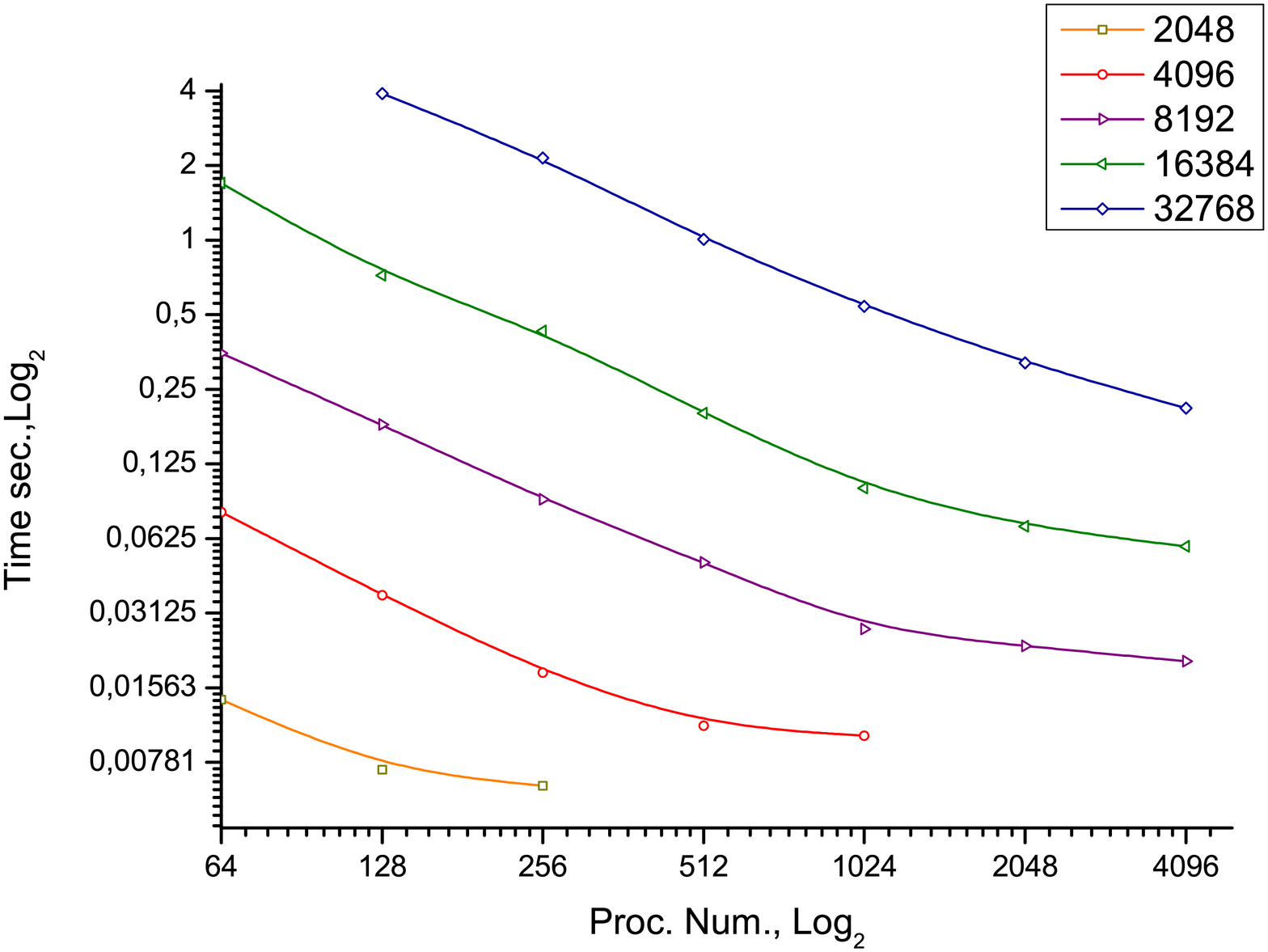}\hfill
\includegraphics[width=0.5\textwidth]{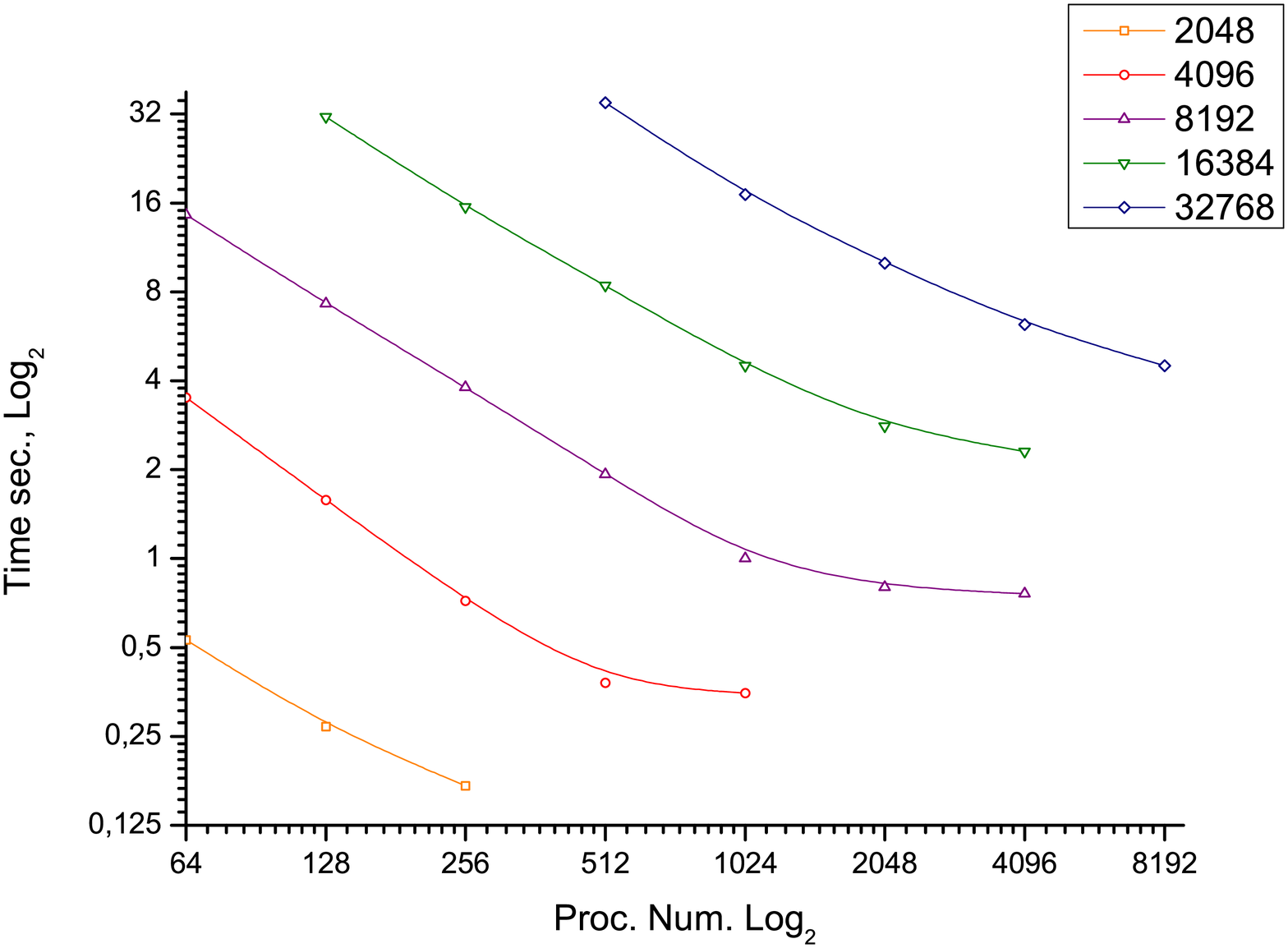}\\
\hfill \parbox[t]{0.4\textwidth}{\small (a) Dependence of
calculation time of one iteration for the conjugate gradient
method for meshes of different resolution.} \hfill
\parbox[t]{0.4\textwidth}{\small(b) Dependence of
calculation time of one cycle for the GMRES$(10)$ method for
meshes of different resolution.}\\
\caption{} \label{TIME}
\end{center}
\end{figure}

Table~\ref{tab1} and Fig.~\ref{TIME}.a represent measurement
results of performance for the conjugate gradient method;
Table~\ref{tab2} and Fig.~\ref{TIME}.b show those for the
GMRES$(10)$ method. Implementing these algorithms, we achieved a
nearly linear dependence of the speedup on the number of
processors for meshes of different resolution. Within one
calculation, we managed to involve a considerable number of
processors (from $1024$ to $8192$) with an efficiency of $90\%$ to
$50\%$, respectively. The achieved performance and scalability are
provided due to using the Dichotomy Algorithm in the context of
the parallel preconditioner inversion. Thus, the algorithm will
substantially increase the efficiency of usage of supercomputer
computational resources in solving elliptic equations. Thus, in
solving applied geophysics problems.

\begin{table}[!h] \center \small
\begin{tabular}{lcccccccccc}
  \hline
  size & \multicolumn{2}{c}{2048x2048} &\multicolumn{2}{c}{4096x4096}& \multicolumn{2}{c}{8192x8192}&\multicolumn{2}{c}{16384x16384}&\multicolumn{2}{c}{32768x32768} \\ \hline
   NP & $ \mathrm{T}$&  $ \mathrm{S}$  &  $ \mathrm{T}$&  $ \mathrm{S}$  & $ \mathrm{T}$&$ \mathrm{S}$&$ \mathrm{T}$&$ \mathrm{S}$&$ \mathrm{T}$&$ \mathrm{S}$\\
  \hline
  64 &1.4e-02 &-&   8.0e-02&     - & 3.5e-01   & -&  1.7      &-&-&-\\
  128 & 7.3e-03&122& 6.3e-02  &138&   1.8e-01   &124&   7.2e-01    &151&3.9&-\\
  256 & 6.3e-03&\underline{142}& 1.8e-02  &284&   9.0e-02   &254&   4.3e-01    &253&2.15&172\\
  512 & -&-& 1.1e-02  & 465 & 5.0e-02   &448&   2.0e-01     &544&1.01&463 \\
  1024 &-&-& 1.0e-02  & \underline{512} & 2.7e-02   &829&       1.0e-01&1088&5.4e-01&924\\
  2048 &- &-&-        &-    & 2.3e-02   &973&7.0e-02&1554&3.2e-01&1560\\
  4096 &- &-&-        &-    & 2.0e-02   &\underline{1120}&5.8e-02&\underline{1875}&2.1e-01&\underline{2377}\\ \hline
\end{tabular}
\caption{Calculation time ($\mathrm{T}$) and speedup
($\mathrm{S}$) versus the number of processors for one iteration
of the CG method.} \label{tab1}
\end{table}

\begin{table}[!h]
\center \small
\begin{tabular}{lcccccccccc}
  \hline
  size & \multicolumn{2}{c}{2048x2048} &\multicolumn{2}{c}{4096x4096}& \multicolumn{2}{c}{8192x8192}&\multicolumn{2}{c}{16384x16384}&\multicolumn{2}{c}{32768x32768} \\ \hline
   NP & $ \mathrm{T}$&  $ \mathrm{S}$  &  $ \mathrm{T}$&  $ \mathrm{S}$  & $ \mathrm{T}$&$ \mathrm{S}$&$ \mathrm{T}$&$ \mathrm{S}$&$ \mathrm{T}$&$ \mathrm{S}$\\
  \hline
  64 &0.53 &-&   3.51&    - &14.6 & -&   -    &-&-&-\\
  128 &0.27&125& 1.58  &142&7.3  &128& 31.3&-&-&-\\
  256 & 0.17&\underline{200}&0.72  &312&  3.8 &245& 15.5&258&-&-\\
  512 & 0.32&106& 0.38 & 591 & 1.93  &484 & 8.4  &476&35& -\\
  1024 &-&-& 0.35  & \underline{641} & 1  &934& 4.5&890&17.1&1047\\
  2048 &- &-&0.5       &450    &0.8    &1168&2.63&1523&9.62&1862\\
  4096 &- &-&-        &-    & 0.76&\underline{1229}&2.3&\underline{1741}&5.7&3132\\
  8192 &- &-&-        &-    &-&-&-&-&4.15&\underline{4318}\\ \hline
\end{tabular}
\caption{Calculation time ($\mathrm{T}$) and speedup
($\mathrm{S}$) versus the number of processors for one cycle of
the GMRES$(10)$ method.} \label{tab2}
\end{table}

Dupros et al. \cite{elast:parallel} considered a parallel
algorithm for solving the dynamic problem of the elasticity
theory. They demonstrated the possibility of using $1024$
processors. Up to $256$ processors, the authors have obtained a
high speedup; however, the algorithm efficiency was lower in the
range from $256$ to $1024$ processors. Our algorithm ensures a
high efficiency in the range from $64$ to $8192$ processors, the
software implementation been much simpler.

As a result of numerical experiments it has been found that the
execution time of the first iteration for the CG and GMRES$(k)$
methods is several times greater than that of subsequent ones.
This is explained by application of dynamic optimization of
interprocessor communications on the level of MPI-Reduce("+")
after repeated execution of the main communication operation "+"\
in the context of the Dichotomy Algorithm. In this case, due to
the associative addition, the order of processor exchanges is set
such that to minimize as much as possible the communication time.
Thus, the possibility of applying the algorithms of dynamic
optimization of the communication interactions ensures the high
performance of the Dichotomy Algorithm. We should note that for
the cyclic reduction method, a fixed order of elimination of
unknowns prevents optimization of the communication interactions
to a full extent. For this reason, in practice, the Dichotomy
Algorithm possesses a much higher performance than the cyclic
reduction method.

It is known that the efficiency of variational methods for solving
SLAEs on a supercomputer decreases because of intensive
communication interactions while computing$\|\cdot\|$ on
distributed data. This problem can be solved by means of
modifications of the known algorithms \cite{Dongarra:NumericalLA}.
Let us compare estimates of the communication time for computing
$\|\cdot\|$   for the CG and GMRES$(k)$ methods and the
communication time of the Dichotomy Algorithm:
    $$
    \begin{array}{l}
    T^{\|\cdot\|,\,all-reduce}_{p}=2\log_2(p)\alpha+\frac{p-1}{p}\left(\gamma+2\beta\right),\\\\
    T^{Dichotomy}_{p}=\alpha\left[\log_2(p)+1\right]\log_2(p)+l\left(\log_{2}(p)-\frac{p-1}{p}\right)\left(\gamma+2\beta\right).
    \end{array}
    $$

From this it follows that for computer systems with a low latency
($\alpha$) and for \mbox{$l\gg 1$}, the communication time for
calculating  $\|\cdot\|$ is insignificant, compared to the
communication costs of the Dichotomy Algorithm. Thus, the chosen
precondition procedure does not need modifications of the CG and
GMRES$(k)$ methods.

\subsection{Acoustic Waves}
The use of mesh methods in spatial derivative approximation cause
a numerical effect called a "phase error"\ \cite{FEM1}. In
modeling wave propagation processes for long instants of time,
this effect determines considerably the accuracy of the solution.
For this reason, from the view of practice, an urgent problem is
choosing the number of mesh nodes per characteristic wavelength.
The high performance of the proposed algorithm allows one to
estimate the accuracy of solution for meshes with a high
resolution ($h_{\alpha}=1/100\lambda\div1/150\lambda$).

Tables 1 and 2 show that calculation of acoustic waves requires
much less count time than that of elastic waves. Hence, we first
consider the problem of modeling acoustic wave propagation in a
homogeneous medium $\rho,\kappa\equiv const$. This made it
possible to investigate the accuracy of solution for meshes with
much more nodes and with less computational costs.

For problem (\ref{acoustic_problem}), a point source was situated
at the origin of coordinates; the time dependence was given as

\begin{equation}
\label{functiont}
 f(t)=\exp\left[-\frac{(2\pi
f_0(t-t_0))^2}{\gamma^2}\right]\sin(2\pi f_0(t-t_0)),
\end{equation}
where $f_0=30\mathrm{Hz},\;t_0=0.2s,\;\gamma=4$.

Approximation of Eq. (\ref{laguerre_h})  was done on the
uniform mesh $\bar{\omega}$ with \mbox{$N_1=N_2=2^k$} nodes,
\mbox{$k=\{12,13,14,15\}$}. The number of addends in series
(\ref{series_lag}) was $n=3000$; the expansion parameters were
$\alpha=9,\;h=400$. The distances were measured in wavelength
$\lambda$ .

\begin{figure}[!htb]
\begin{center}
\includegraphics[width=0.45\textwidth]{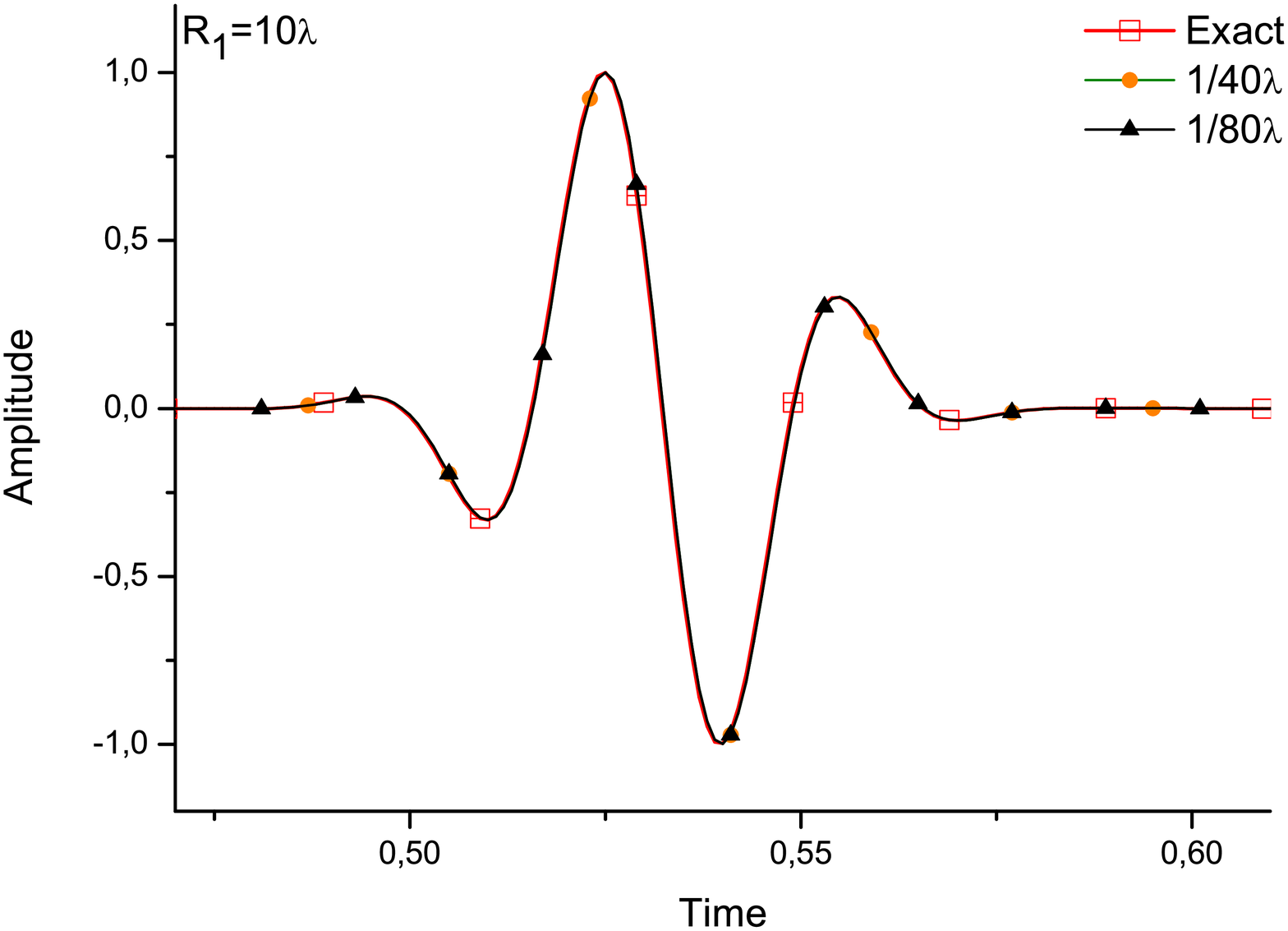}\hfill
\includegraphics[width=0.45\textwidth]{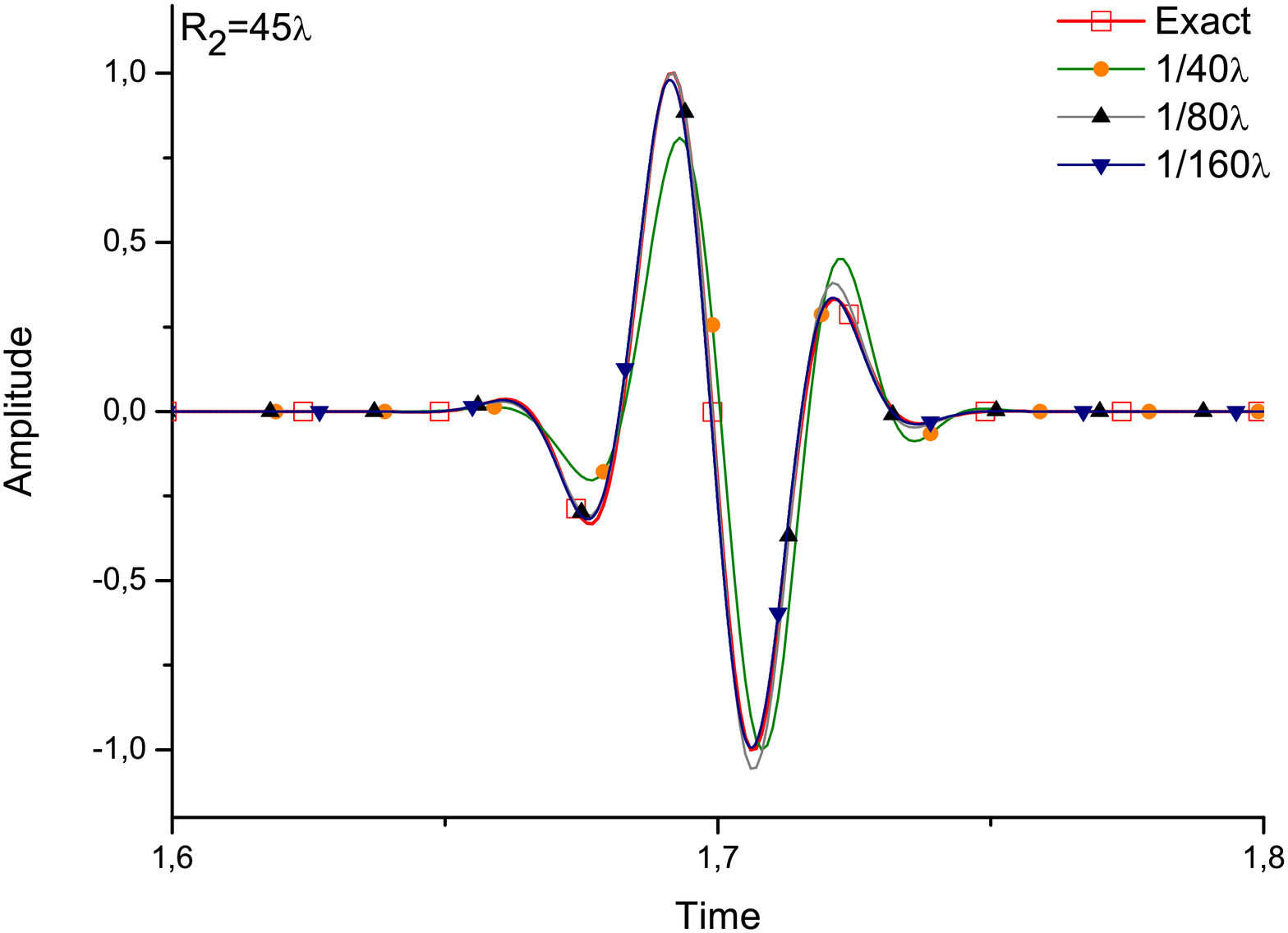}\\
\parbox[t]{0.5\textwidth}{\center a)}\hfill
\parbox[t]{0.5\textwidth}{\center        b)}\\
\includegraphics[width=0.45\textwidth]{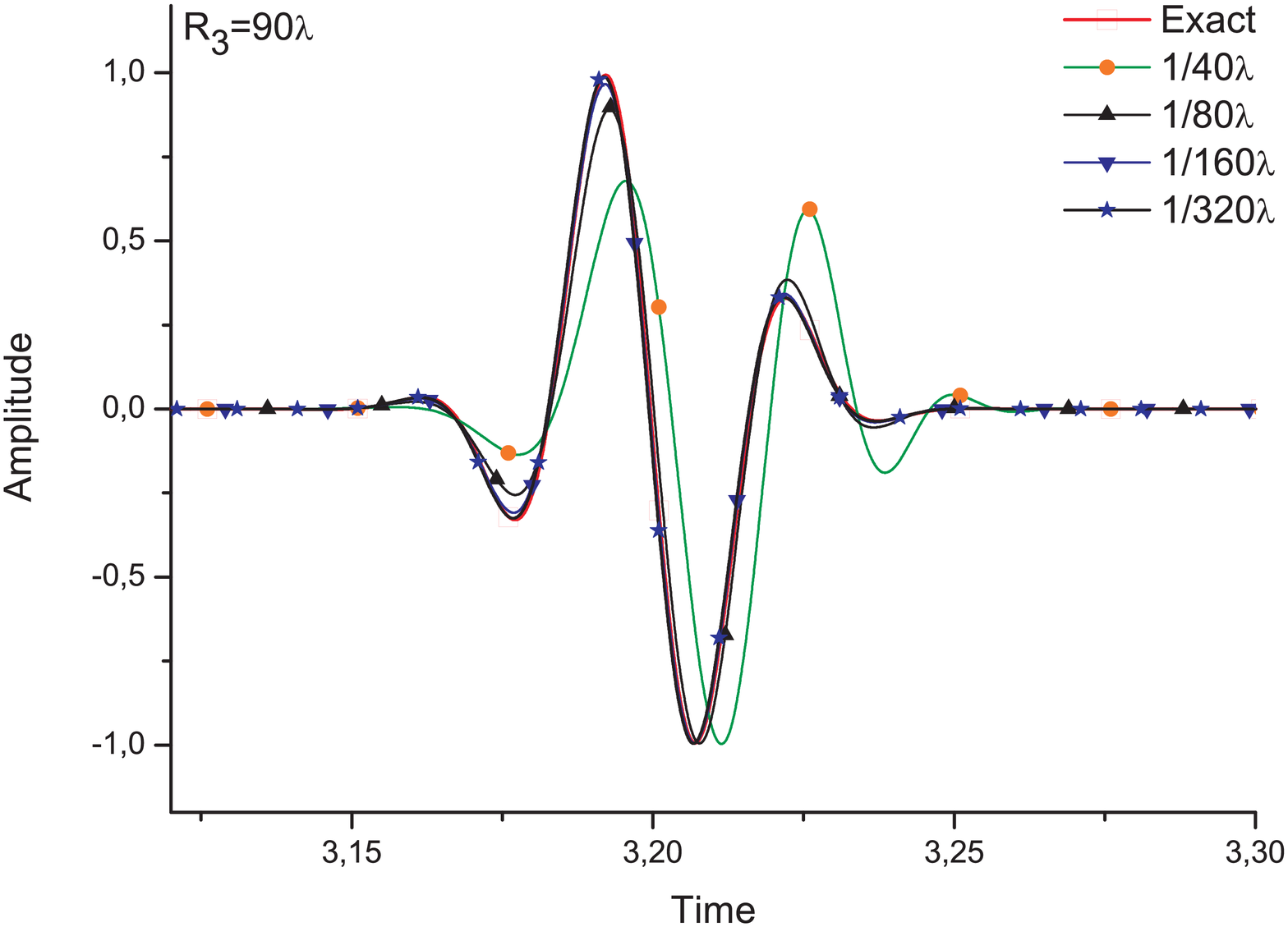}\hfill
\includegraphics[width=0.45\textwidth]{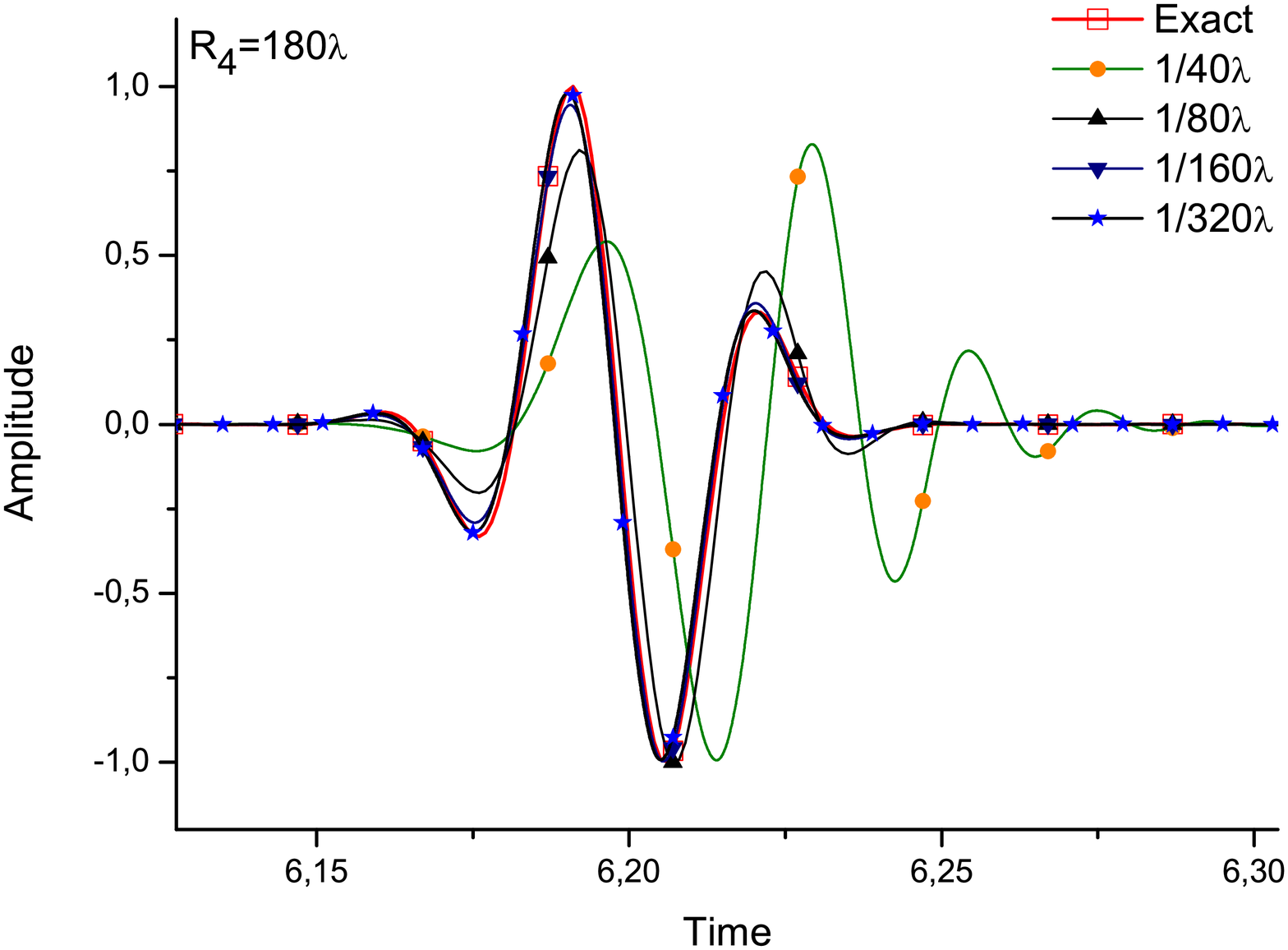}\\
\parbox[t]{0.5\textwidth}{\center c)}\hfill
\parbox[t]{0.5\textwidth}{\center        d)}\\
\caption{ The time dependence of solution $u({\bf x}_i,t)$ for the
acoustic equation, where \mbox{${\bf x}_i=(\mathrm{R}_i,0)$},
\mbox{$i=1,2,3,4$}.} \label{main_pic223}
\end{center}
\end{figure}

The time dependencies of the wavefield amplitudes for four
receivers situated on the free surface at different lengths from
the source are represented in Fig.~\ref{main_pic223}. It is seen
that the accuracy of the obtained solution for different instants
of time depends substantially on the number of mesh nodes per
characteristic wavelength. For instance, for first instants of
time, for achieving a reasonable calculation accuracy, the mesh
with a space step $h_r=h_z=1/40\lambda$ is sufficient
(Fig.~\ref{main_pic223}.a); for longer time intervals it is
required to decrease the mesh step in order to keep a reasonable
level of the calculation accuracy.

\begin{figure}[!htb]
\begin{center}
\includegraphics[width=0.5\textwidth]{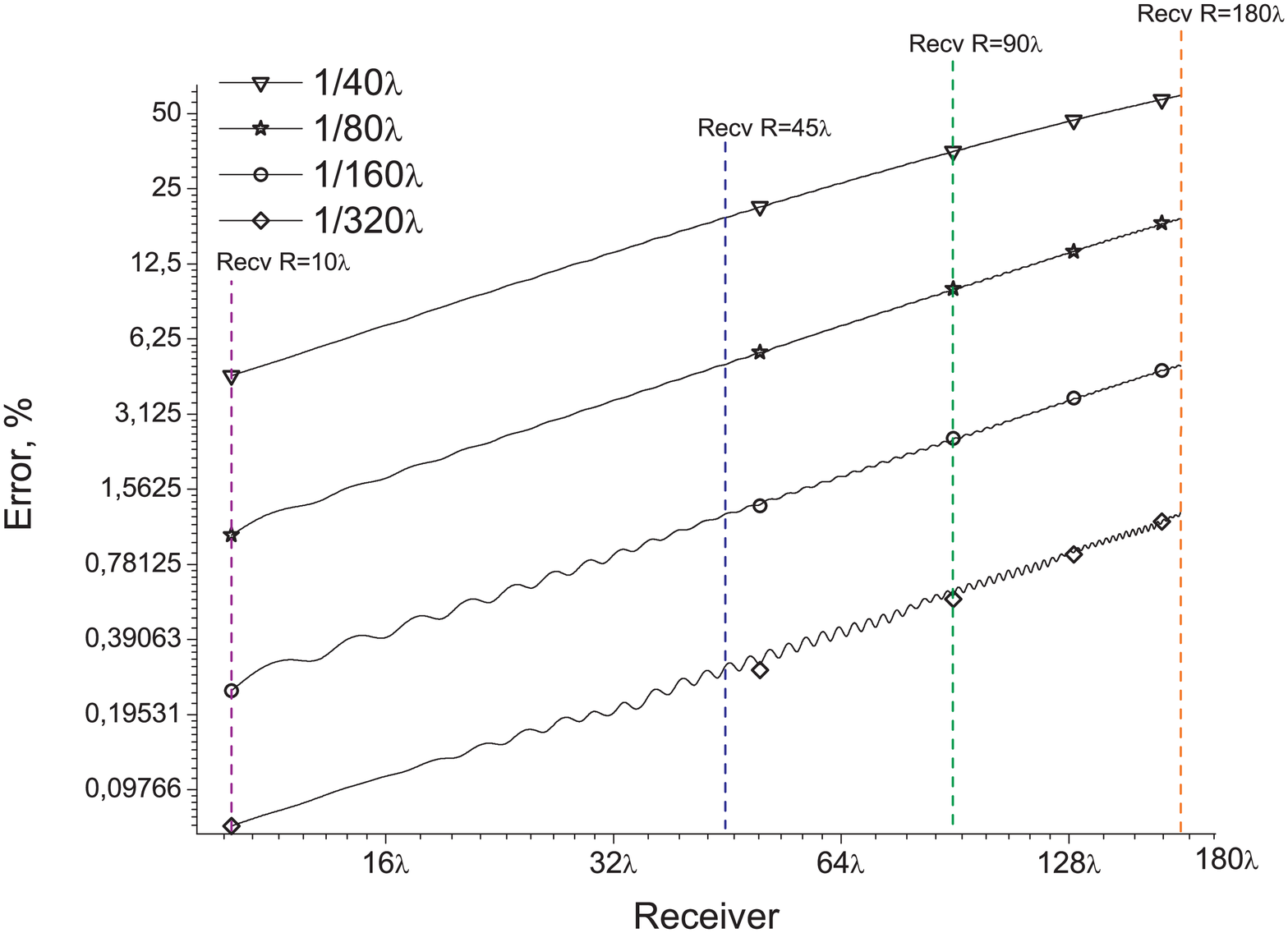}
\end{center}
\caption{Dependence of the solution accuracy on the position of
receiver (\ref{accuracy}) for the meshes of different resolution.}
\label{pic:phaseerror}
\end{figure}

Figure~\ref{pic:phaseerror} represents dependence of the accuracy
of the obtained solution on the receiver position for meshes with
different resolution:
$$
 \epsilon({\bf
x}_i)=\sqrt{\frac{\int_{0}^{t_1}\left[u_{exact}({\bf
x}_i,t)-u_h({\bf x}_i,t)\right]^2\mathrm{dt}}{\int_{0}^{t_1}
\left[u_{exact}({\bf x}_i,t)\right]^2\mathrm{dt}}},\; {\bf
x}_i=(ih_r,0),\; i=1,..,N_1 \label{accuracy},
$$

where $\displaystyle
u_{exact}(r,0,t)=\frac{1}{2\pi}\frac{f(t-r/\sqrt{\kappa/\rho})}{r}$
is the exact solution and $u_h$ is the numerical solution obtained
on the mesh with the space step $h_r=h_z=h$.

It can be easily found that when the time interval of modeling is
increased $m$ times, the space step must be decreased
$\approx\sqrt{m}$ times; this agrees with theoretical estimates
for approximation methods of the second order of
accuracy\cite{FEM1}. Thus, for acoustic wave simulation for long
time intervals, it is necessary to use meshes with a sufficient
number of nodes in order that numerical effects caused by the
model resolution do not predominate.

We will note that for solving the problem, higher-order schemes
\cite{Kim1,Approx:Mech1} are suitable. In the context of the
parallel algorithm, increasing approximation order does not cause
loss of efficiency because the preconditioning for higher-order
schemes can be done with the second order. Naturally, in this case
the number of iterations for the CG and GMRES$(k)$ methods for
achieving the desired accuracy will be a bit more, but the
behavior of the dependence of the speedup on the number of
processors will not change.

\subsection{Solid layer over Solid Half Space}

Although early results on elastic wavefield modeling have been
obtained long ago \cite{Kelly:history}, \cite{Kelly:Love},
however, they were rather qualitative because of a large step of
the space mesh $h=1/5\lambda\div1/2\lambda$. Considerably
increased computer performance and also development of
multiprocessor computer systems have made it possible to increase
the calculation accuracy \cite{Bruno1,Ito1,MPI_elastic}. However,
in spite of available theoretical estimates of the dependence of
solution accuracy on mesh step \cite{FEM1}, the problem of
practical choosing a space step of meshes is still urgent. By
solving the acoustic equation, it was illustrated that
calculations for long instants of time require meshes with many
nodes. Taking into account that the proposed parallel algorithm
possesses high performance, we will analyze issues of accuracy for
problem (\ref{elastic_problem}) for meshes
$h=h_r=h_z=\{1/10\lambda_s,1/20\lambda_s,1/45\lambda_s,1/90\lambda_s\}$,
where  $\lambda_{s}=\min V_s/f_0$. Here, $V_s$ is the $S$-wave
propagation velocity and $f_0$ is the source frequency.

\begin{figure}[!h]
\begin{center}
\includegraphics[width=\textwidth]{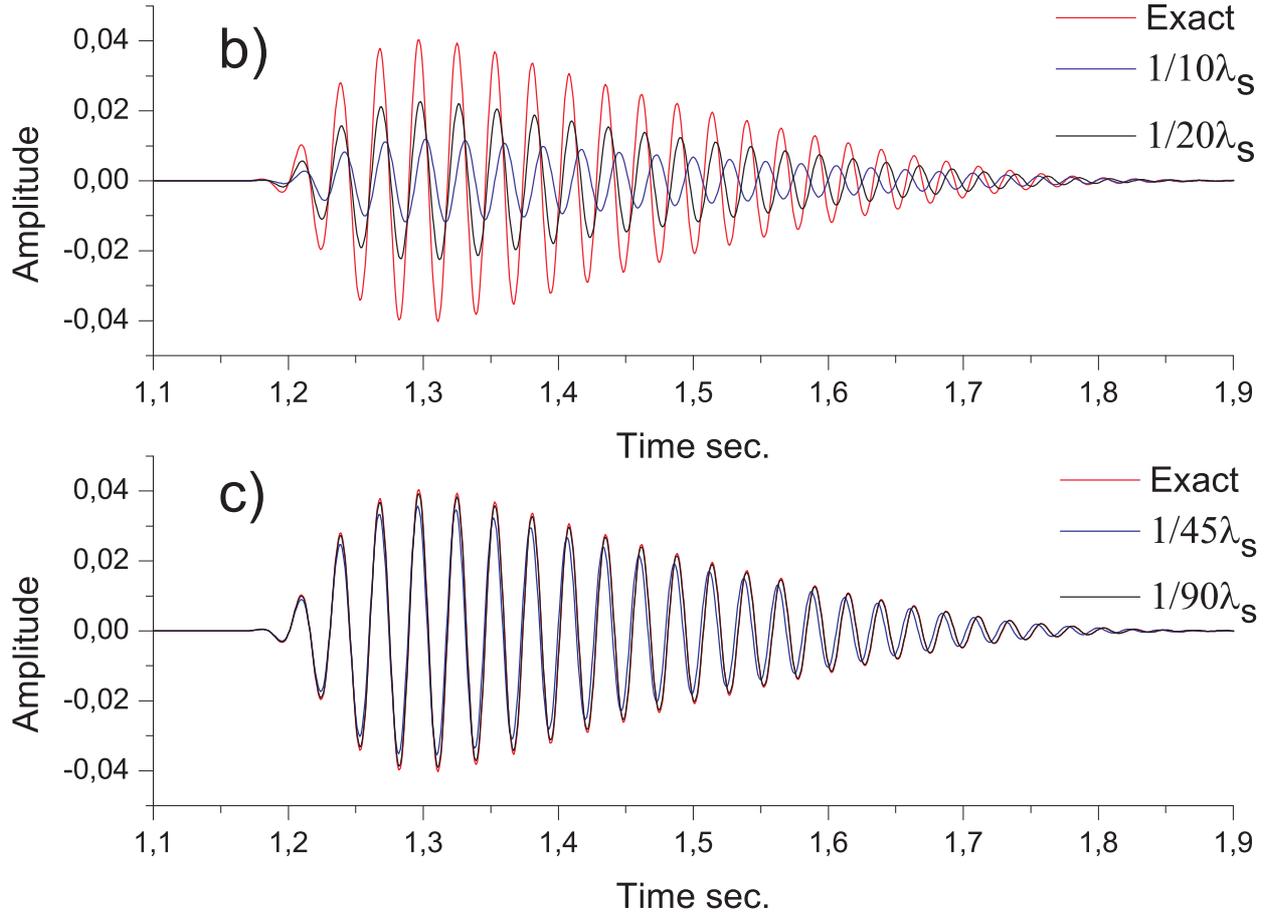}
\caption{(a) A snapshot for the displacement vector component
$u_z(z,r)$ at $t=3s$ in the presence of a thin layer. The time
dependence of amplitude  $t\in [1.1,1.9]s$ for detector receiver
$u_z(1500m,0)$ for calculations on the meshes: (b) $N_r\times
N_z=\{4096\times4096,8192\times8192\}$; (c) $N_r\times
N_z=\{16384\times16384,32768\times32768\}$.} \label{main_pic233}
\end{center}

\end{figure}

Let us consider a problem on elastic wave propagation in a thin
layer whose seam thickness is comparable with the wavelength
(Fig.~\ref{main_pic233}.a). The wavefield source is a source of
the type of "center of pressure" \cite{Achenbach}:

\begin{equation}
F_r=\frac{1}{2\pi}\frac{\mathrm{d}}{\mathrm{d}r}\left[\frac{\delta(r)}{r}\right]\delta(z-d),\quad
F_z=\frac{1}{2\pi}\frac{\delta(r)}{r}\frac{\mathrm{d}}{\mathrm{d}z}\delta(z-d).\label{centre_preassure}
\end{equation}
The time dependence of the pulse $f(t)$ was determined in
(\ref{functiont}), where $f_0=30\mathrm{Hz},\;t_0=0.2s$ and
$\gamma=4$. The source is placed at the depth $d=10m$. In the
calculations, the instants of time were $t\in(0,5] $s. The number
of terms of series (\ref{series_lag2}) was $n=2000$ with the
parameters $\alpha=8$ and $h=600$.

In problems of simulation of wavefields, in particular, seismic
ones, the governing factor is choosing model problems to estimate
accuracy of numerical algorithms. A common method applied for
layered media is the method proposed in \cite{3f,4f} and extended
in \cite{1f,2f}, etc. The drawback of the method is that it
introduces interference (artifacts). The use of minor matrices
\cite{5f} made it possible to extend applicability of the matrix
method, but did not eliminated artifacts. In the present paper, we
evaluate the accuracy of the proposed parallel algorithm with the
use of the approach described in \cite{6f,7f}. The essence of the
method is that the sought-for boundary-value problem of second
order in the spectral domain is reduced to two Cauchy problems of
first order, to which there exists a stable analytical solution.
As a consequence, this made it possible to remove all constraints
on the powers of layers, frequencies, and recording systems.
Comparison of the modeling results and results obtained by means
of the analytical method made it possible to estimate the
dependence of the numerical solution accuracy on the space mesh
step.

\begin{figure}[!h]
\begin{center}
\includegraphics[width=\textwidth]{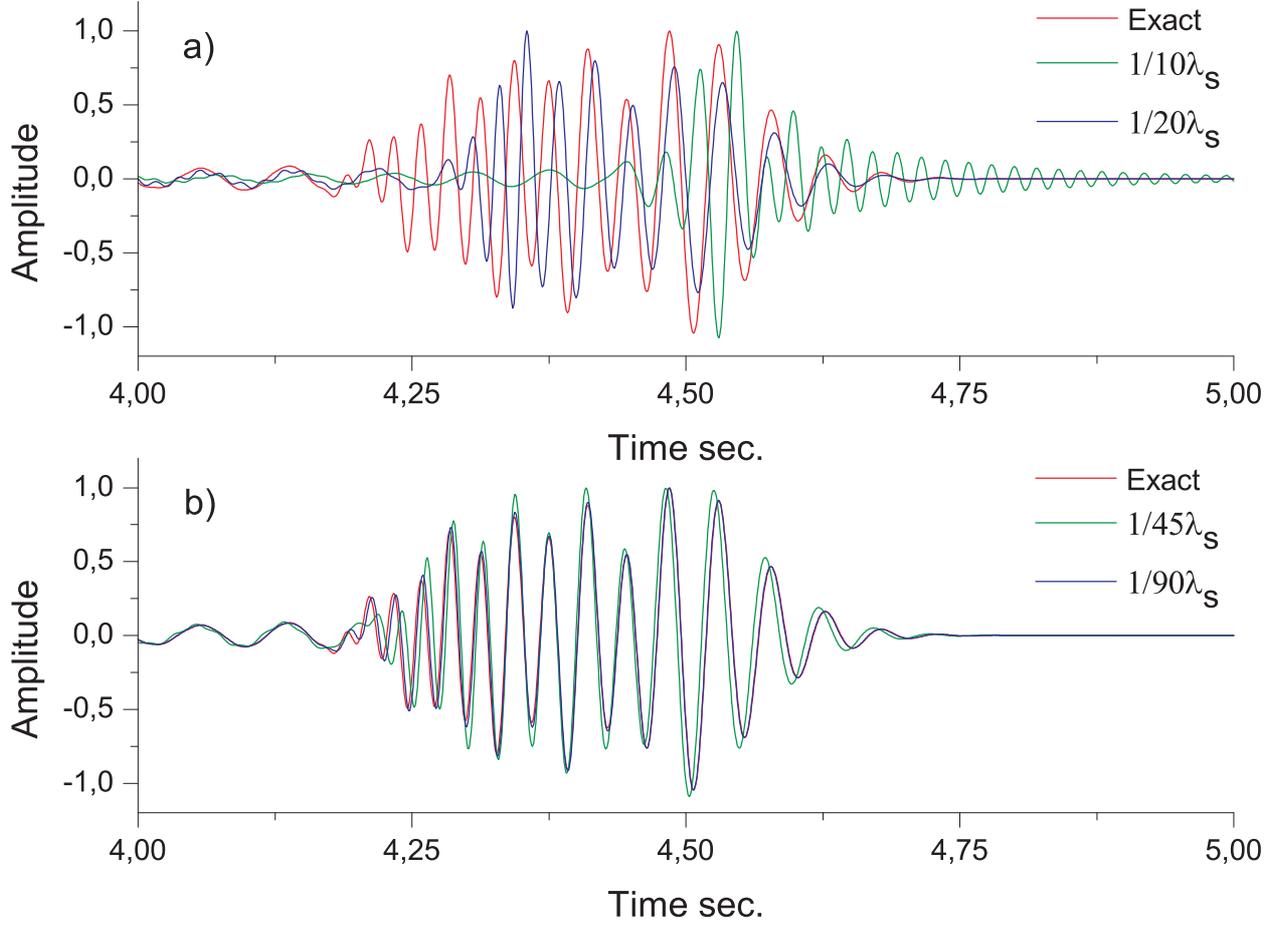}
\caption{The time dependence of amplitude $t\in [4,5]s$ for
detector receiver $u_z(1500m,0)$ for calculations on the meshes:
(a) $N_r\times N_z=\{4096\times4096,8192\times8192\}$; (b)
$N_r\times
N_z=\{16384\times16384,32768\times32768\}$.}\label{main_pic234}
\end{center}
\end{figure}

The medium model and a snapshot of the wavefield for the component
$u_z(r,z)$  at $t=3$s. are represented in
Fig.~\ref{main_pic233}.a. Figures~\ref{main_pic233}.b,c and
Figs.~\ref{main_pic234}.a,b show the component $u_z$ as a function
of time for a receiver situated on the free surface at $r=1500$.

Figure~\ref{main_pic233}.b and Fig.~\ref{main_pic234}.a evidence
that the meshes with steps $1/10\lambda_s$ and $1/20\lambda_s$
ensure no accuracy. Thus, the mesh with $1/45\lambda_s$  or
$1/90\lambda_s$ can ensure an acceptable level of accuracy for
initial instants of times (Fig.~\ref{main_pic233}.b and
Fig.~\ref{main_pic234}.a). For final instants of time
(Fig.~\ref{main_pic233}.c and Fig.~\ref{main_pic234}.b) the
calculations have to be done on the mesh with $1/90\lambda_s$ .
The results of evaluating the accuracy of solution for the problem
of the elasticity theory are in good agreement with the results
obtained for the acoustic equation. Therefore, we can state that
the main factor determining the accuracy of solution in wave
process simulation is the number of mesh nodes per wavelength.
Also, modeling of real space-time scales requires modeling meshes
with a higher resolution.
\subsection{Marmousi}

\begin{figure}[!h]
\begin{center}
\includegraphics[width=1\textwidth]{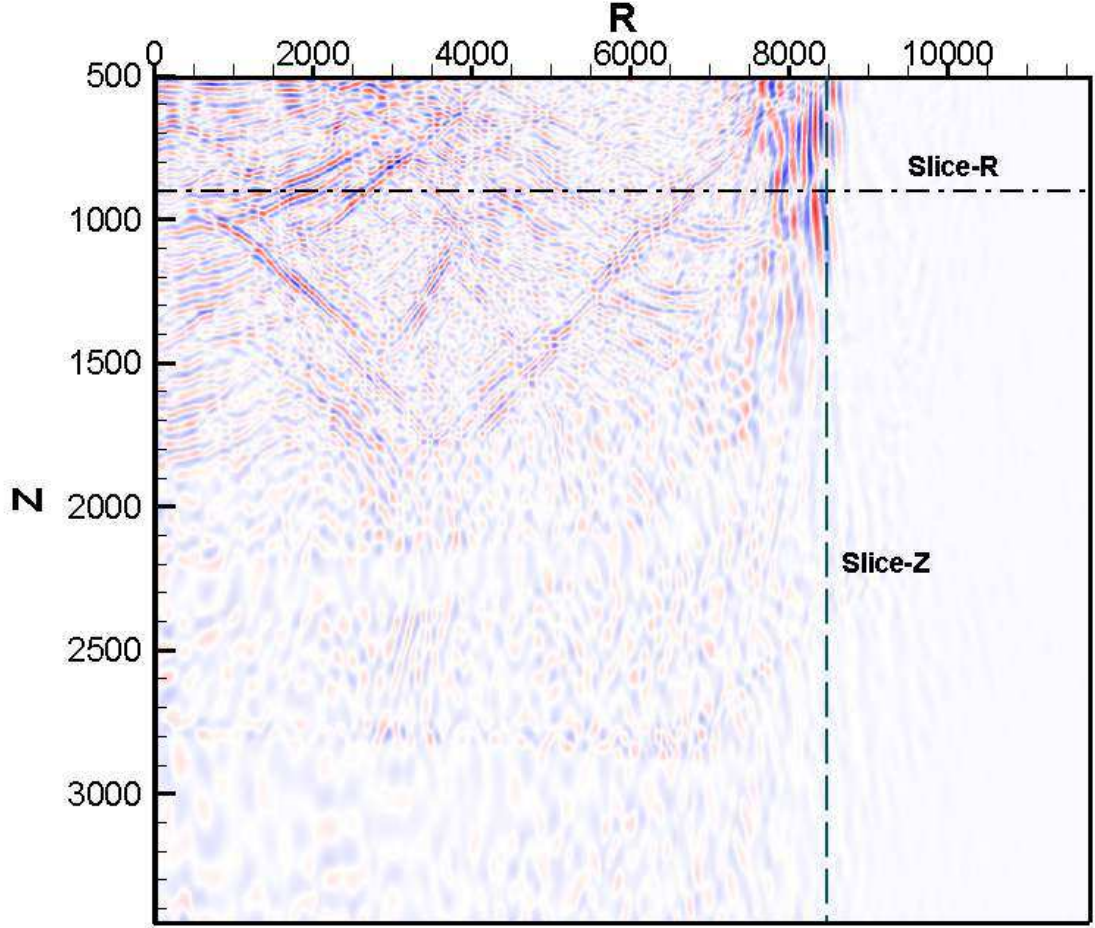}
\caption{(a) Marmousi model ($z\geq 500m$): P,S -- wave
velocities; (b) a snapshot for the displacement vector component
$u_r(r,z)$ at $t=6$ s.} \label{pic:marmslice}
\end{center}
\end{figure}

For illustrating the ability of the proposed algorithm to perform
(for an acceptable time) elastic wavefield simulation for real
application problems, we will consider problem
(\ref{elastic_problem}) for the Marmousi medium
(Fig.~\ref{pic:marmslice}.a) \cite{Marmousi2}.

\begin{figure}[!h]
\begin{center}
\includegraphics[width=0.8\textwidth]{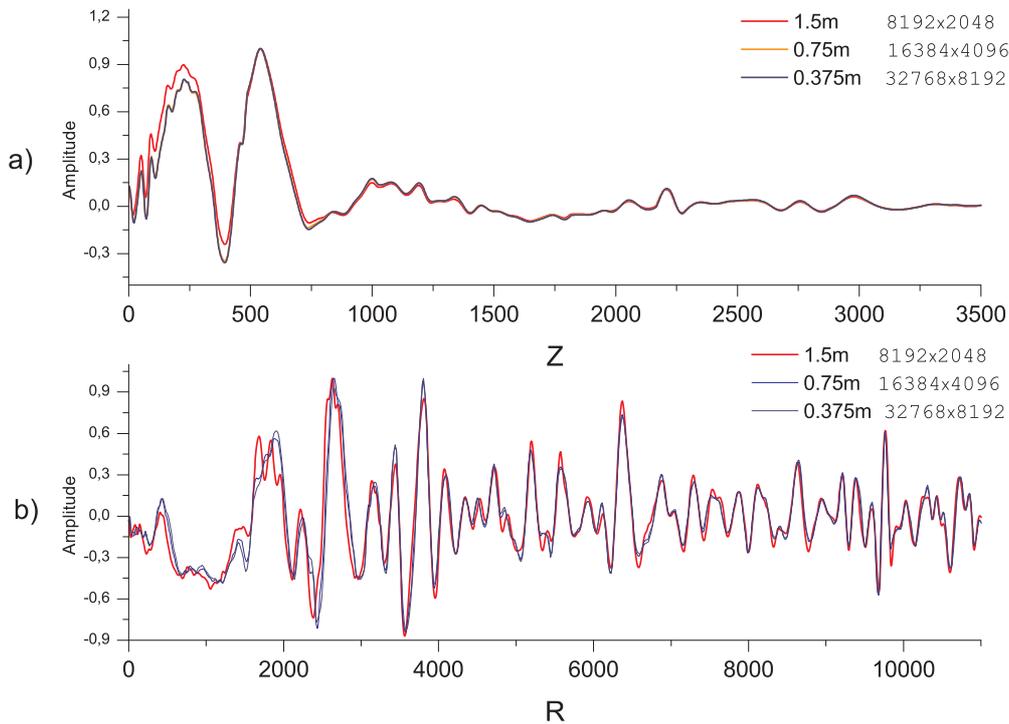}
\caption{Dependence of the field amplitude for the component $u_r$
at $t=6$ along straight lines (a) Slice-Z and (b) Slice-R.}
\label{pic:marmslice2}
\end{center}
\end{figure}

The calculations were done for $t \in(0,6] $s on the meshes with
$N_{r}\times N_z=\{8192 \times 2048,16384 \times 4096,32768 \times
8192\}$ nodes, which corresponded to a space step
$h_r=h_z=\{1.5m,0.75m,0.375m\}$. The wavefield was modeled from a
source of the type of center of pressure
(\ref{functiont}),(\ref{centre_preassure}) with the parameters
$f_0=10\mathrm{Hz},\;t_0=1s$ and $\gamma=4$. The number of terms
of series (\ref{series_lag2}) was $n=1200$ with the parameters
$\alpha=8$ and $h=300$. Figure~\ref{pic:marmslice}.b shows a
snapshot of a wavefield for the displacement vector component
$u_r(r,z)$ at $t=6$s; Figs.~\ref{pic:marmslice2}.a,b represent
dependencies of the component $u_r(r,z)$ along straight lines
"Slice-R"\ and "Slice-Z"\ . Comparing results obtained for
different meshes, we conclude that an acceptable accuracy for the
final instant of time is achieved in calculations with the mesh
space step $h_r=h_z=0.75$m, which corresponds to $\{N_r\times
N_z\}=\{16384 \times 4096 \}$ nodes. For this mesh, according to
Table~\ref{tab3}, computing will take $4.7$ hours with $1024$
processors. The efficiency is about $90\%$. We should note that
sometimes it is reasonable to perform computing using more
processors, but with a lower efficiency, because in this case the
amount of storage is increased. This makes it possible to choose
greater values of $k$ for the GMRES$(k)$ method and, thereby
ensure a higher convergence rate\cite{Saad}.

\begin{table}[!h]
\center \small
\begin{tabular}{l|ccccc}
  \hline
   NP   &$256$&  $512$  &  $1024$&  $2048$  & $ 4096$\\
   \hline
   Mesh & &&Time Hours\\
  \hline
  8192$\times$ 2048 & 3.1  & 1.6  &   1.4 & 2.4  & 4.2\\
  16384$\times$4096 &17 &8.6 &4.7 &3.8 &4\\
  32768 $\times$ 8192 &68 &34 & 19.6&12 &11\\ \hline \end{tabular}
  \caption{Calculation time versus the number of processors for meshes of different resolution for the Marmousi medium.}
  \label{tab3}
\end{table}

The numerical experiments have shown that the proposed algorithm
enables not only to perform modeling, but also to solve real
application geophysics problems. For doing so, both a
supercomputer with a moderate number of processors ($64\div256$)
and multiprocessor systems integrating thousands of computing
elements are used effectively.

\section{Conclusions}

We have proposed the parallel algorithm for solving an acoustic
equation and dynamic problem of the elasticity theory in a
cylindrical coordinate system (2.5 D). The Laguerre time transform
was used to perform changing from the initially boundary-value
problem to the problem of inversion of the same elliptic
second-order operator for different right-hand sides. The
difference equations resulting from elliptic operator
approximation were solved by the CG method or the GMRES method.
Choosing the Laplace operator as the preconditioning one allowed
for a high convergence rate of the iterative process for media
with a moderate contrast.

The nearly linear dependence of the speedup and the high
scalability of the parallel algorithm on the number of processors
were ensured due to the Dichotomy Algorithm in the context of the
variable separation method for inverting the preconditional
operator. The proposed algorithm has validated its efficiency in
calculations with $64$ to $8192$ processors. Thus, the high
performance of the Dichotomy Algorithm and its simple
implementation enable efficient parallelization of economic
numerical procedures that require multiple solution of tridiagonal
systems of equations.

The main conclusion is that the wave process modeling for longer
time intervals requires increasing number of space mesh nodes.
This causes the necessity of applying high-performance computer
systems for solving application problems. It has been shown that
the proposed parallel algorithm, which based on the known economic
numerical methods and the Dichotomy Algorithm, makes it possible
to efficiently involve thousands of processors within one
calculation. This enables to perform practical calculations for
real models of media, times, and distances with the desired
accuracy.
\newpage
\bibliography{base}
\end{document}